\documentclass[a4paper,12pt]{amsart}
\usepackage{amssymb,amsthm,bbm}
\usepackage{a4wide}
\usepackage{amssymb,amsthm}
\usepackage{graphicx}
\newcommand{\Png}[2]{
\centerline{\includegraphics[height=#1]{#2}}
}
\newcommand{\un}{\,\, \vrule height 10pt width 1pt \,\,}
\usepackage{xypic}
\def\k{\mathbbm{k}}

\def\Q{\mathbbm{Q}}
\def\Z{\mathbbm{Z}}

\def\K{\mathbbm{K}}
\def\L{\mathbbm{L}}
\def\End{\mathrm{End}}
\def\Ker{\mathrm{Ker}\,}

\def\B{\mathfrak{B}}
\def\ST{\widetilde{\mathfrak{S}}}
\def\WB{\widehat{\mathfrak{B}}}
\def\wp{\widehat{\varphi}}
\def\Bur{R_{\mathrm{bur}}}
\def\bur{\rho_{\mathrm{bur}}}
\def\Sym{R_{\mathrm{sym}}}
\def\sym{\rho_{\mathrm{sym}}}
\def\Id{\mathrm{Id}}
\def\SS{\mathfrak{S}}
\def\ss{\sigma}
\newtheorem{prop}{Proposition}
\newtheorem{theo}{Theorem}
\newtheorem{lem}{Lemma}
\newtheorem*{ccor}{Corollary}
\newtheorem*{dfn}{Definition}
\theoremstyle{remark}

\newtheorem*{rems}{Remarks}

\title{Cabling Burau representation}
\author{Christian Blanchet and Ivan Marin}
\date{23st October 2006}
\begin{document}
\maketitle
\begin{abstract}
The Burau representation enables to define many other representations
of the braid group $B_n$ by the topological operation of ``cabling braids''. We show here that
these representations split into copies of the Burau
representation itself and of a representation of $B_n/(P_n,P_n)$.
In particular, we show that there is no gain in terms of
faithfulness by cabling the Burau representation.
\end{abstract}

\section*{Introduction}

The Burau representation is the oldest and most natural
non-trivial representation of the braid group $B_n$. Debates over its
faithfulness were a central concern in the past decades, until
it has been shown to be non faithful for $n \geq 5$, the case $n=4$
remaining open.

A natural question, brought to us by T. Fiedler and S. Orevkov,
is whether it is possible to reduce the size of the kernel
of the Burau representation by the operation of ``cabling braids''.
Here we answer this question for the most natural cabling,
sometimes called parallel cabling.

The answer is negative. More precisely, letting
$\Delta_{n,r} : B_n \to B_{nr}$ denote the morphism of cabling
where each strand is replaced by $r$ parallel strands,
as in the diagram below,

\begin{center}
\resizebox{5in}{!}{\includegraphics{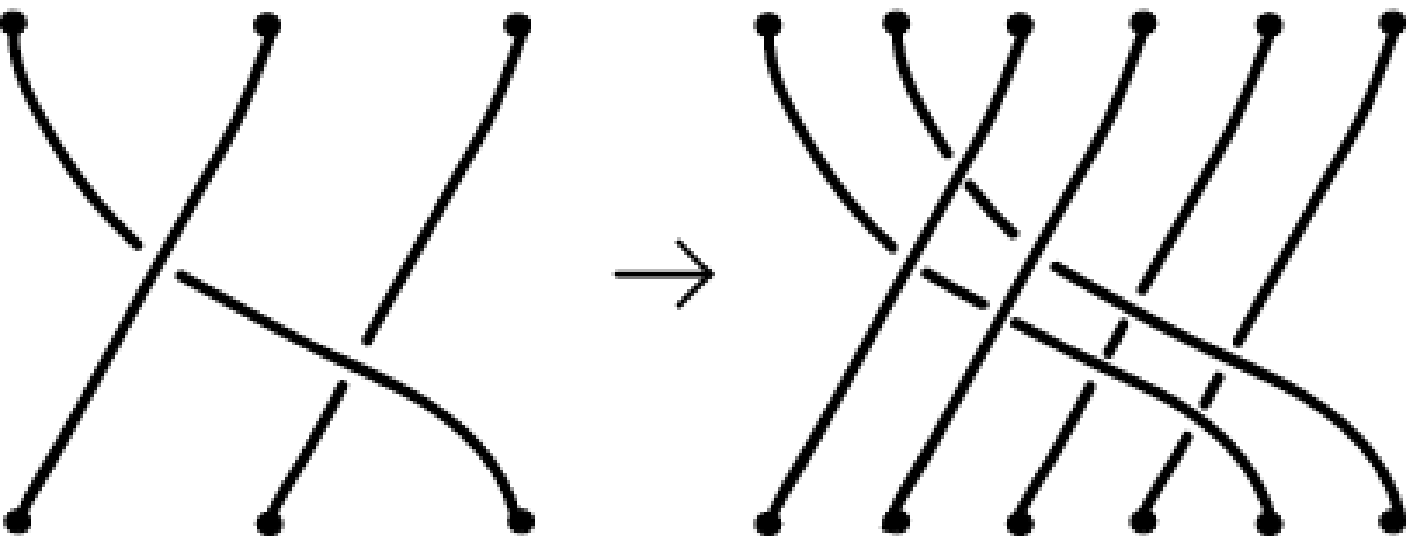}}
\end{center}

\noindent and letting $\Bur : B_{nr} \to GL(V)$
denote the Burau representation of $B_{nr}$, we decompose the representation
$\Bur \circ \Delta_{nr}$ of $B_n$. Our main theorem states that this representation
is semisimple and splits into a variant of the Burau representation
of $B_{nr}$ and copies of an irreducible representation
of $B_n/(P_n,P_n)$, where $P_n$ denotes the pure braid group.
As a consequence, we show that the kernel of $\Bur \circ \Delta_{nr}$
coincides with the kernel of the Burau representation of $B_n$.

Our proof has two ingredients, both based on the use of
Drinfeld associators. These associators define
morphisms from the (group algebra of) the braid group $B_n$
to some ``infinitesimal braids algebra'' $\B_n$.
In the first part we recall the classical extension
of this construction, originally due to Kontsevich, 
in terms of monoidal categories. In particular, we insist
on the compatibility with cabling, for which we give a self-contained
proof.

The second part is based on a correspondance
between representations of $B_n$ and representations of its
infinitesimal couterpart $\B_n$,
for which one can find a detailed account in \cite{REPTHEO}.
We use the fact that the Burau representation corresponds to a very simple
representation of $\B_n$, and that cabling on the
infinitesimal braids is also an easy (additive) operation, to
get the decomposition that we are
interested in.

\section{Drinfeld-Kontsevich functor for parenthesized braids} 
\subsection{Braids and chord diagrams}
The $n$ strands braid group $B_n$ is the fundamental group of the configuration space $Y_n$ of
$n$ points in the complex plane.  It may be convenient to vary the base point, and to consider the fundamental groupoid; parenthesized braids amount to consider
limit configurations in a convenient compactification of $Y_n$.

A parenthesized braid \cite{BNGT} is a braid whose ends are points on the real line, together
with a parenthesization of its bottom end (the domain) and its top end (the range).
 We obtain a groupoid $\mathbf{PaB}$ which is a subcategory in the parenthesized tangle category also called the $q$-tangle category \cite{LM1,LM2}.

Let $\mathcal{A}_n=\mathcal{A}_n(\k )$ be the algebra
of  chords diagrams for $n$-strands
pure braids, over the  scalar field $\k$ of characteristic $0$. As a unital $\k $-algebra $\mathcal{A}_n$ is generated by the $t_{ij}=t_{ji}$, $0\leq i<j\leq n$, represented by a chord between strands numbered $i$ and $j$, with relations (infinitesimal braid relations):
$$[t_{ij},t_{kl}]=0 \text{ if the $4$ indices $i,j,k,l$ are distinct,}$$
$$[t_{jk},t_{ij}+t_{ik}]=0 \text{ if the $3$ indices $i,j,k,$ are distinct.}$$
 The number of chords provides a grading on the algebra $\mathcal{A}_n$, and we denote by 
$\widehat{\mathcal{A}}_n$ the completion with respect to this grading.

We use the notation $(\sigma,D)\mapsto {}^\sigma\! D$ for the natural left action of the symmetric group $\mathfrak{S}_n$  on $\mathcal{A}_n$ (resp. $\widehat{\mathcal{A}}_n$);
 for a generator $t_{ij}$, we have ${}^\sigma(t_{ij})= t_{\sigma(i)\sigma(j)}$. 
 The algebra $\mathfrak{B}_n$ (resp. $\widehat{\mathfrak{B}}_n$) is defined
as the crossed product of $\mathcal{A}_n$ (resp. $\widehat{\mathcal{A}}_n$)
with the symmetric group $\mathfrak{S}_n$. As a module, $\mathfrak{B}_n$ (resp. $\widehat{\mathfrak{B}}_n$) is a free $\mathcal{A}_n$-module (resp. $\widehat{\mathcal{A}}_n$)-module), with basis $\mathfrak{S}_n$.
 In Vassiliev finite type invariants theory, $\mathfrak{B}_n$ is the algebra
of chords diagram for braids. Its completion $\widehat{\mathfrak{B}}_n$
is the natural target for the Drinfeld-Kontsevich functor which we want to
consider now. We will need extra structures on our categories, namely cabling operations and strand removal operations \cite{BNGT}.

\begin{figure}
 \centerline{
\Png{3cm}{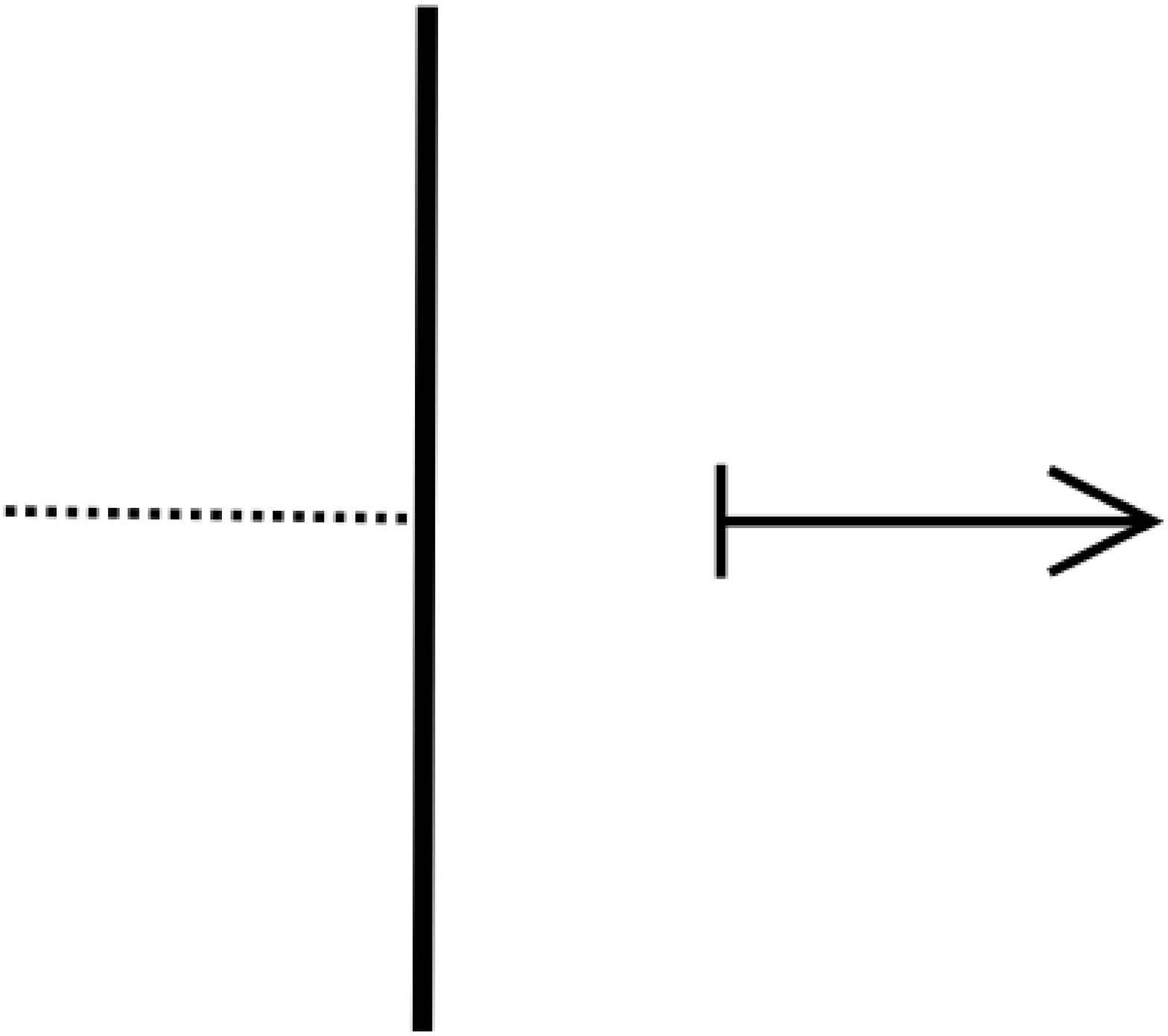}
}
\caption{\label{cablingChord}}
 \end{figure}

\begin{figure}
 \centerline{
\Png{5cm}{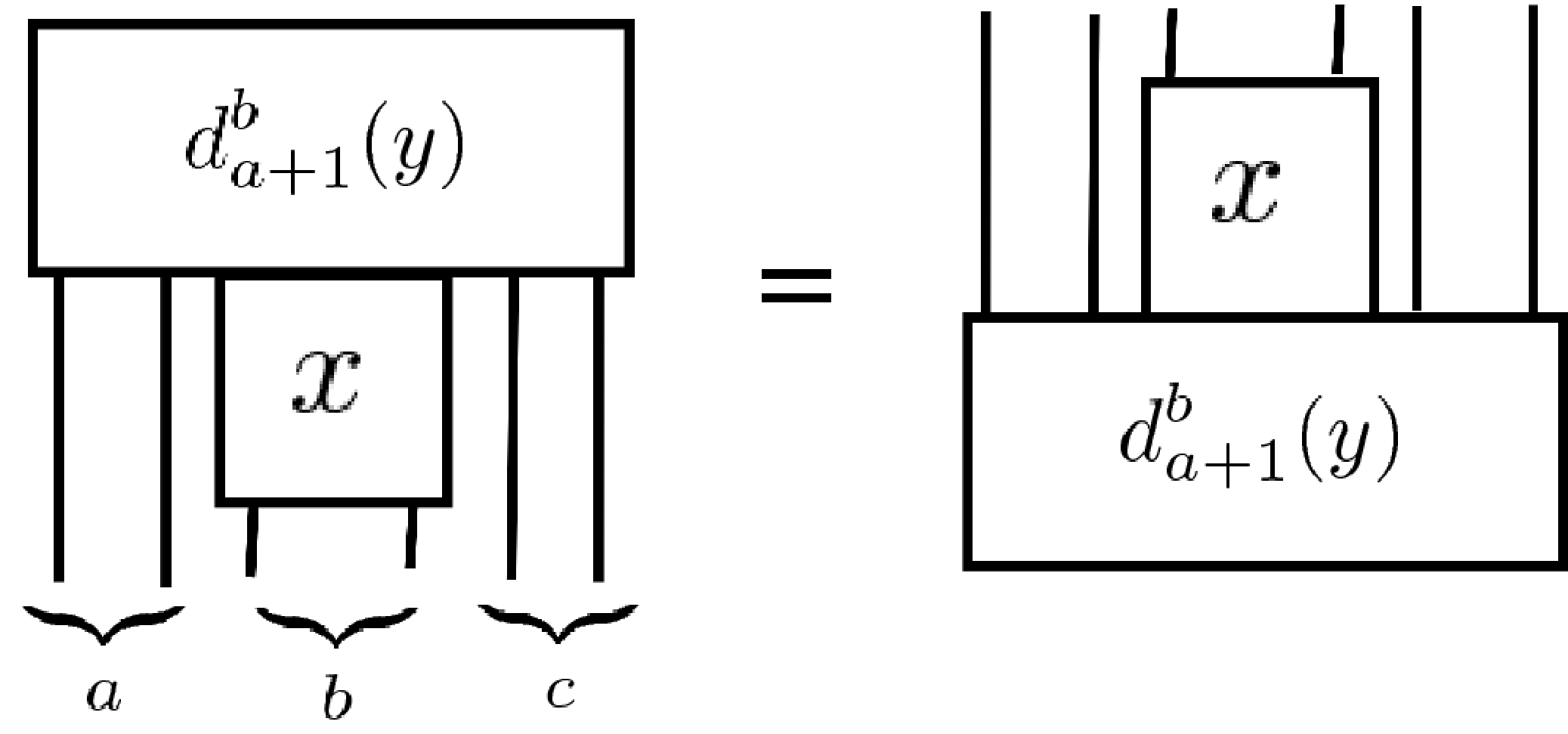}
}
\caption{\label{naturalCabling}}
 \end{figure}

\vspace{5pt}
\noindent {\em Cabling braids:} For a parenthesized braid $B$, $d_i(B)$ is the parenthesized braid obtained from $B$ by doubling the $i$th strand (counting at the bottom).\\[5pt]
\noindent {\em Cabling chord diagram:} For a chord diagram $D$, $d_i(D)$ is the chord diagram obtained form $D$ by doubling the $i$th strand (counting at the bottom). Each chord incident to the
$i$th strand is expanded as depicted in figure \ref{cablingChord}. Note that $d_i\circ d_i=d_{i+1}\circ d_i$.\\[5pt]
\noindent{\em Strand removal on braids:} For a parenthesized braid $B$, $s_i(B)$ is obtained by removing the
$i$th strand (counting at the bottom).\\[5pt]
\noindent {\em Strand removal on chord diagrams:} For a chord diagram $D$, $s_i(B)$ is obtained by removing the
$i$th strand (counting at the bottom) if no chord is incident to this $i$th strand, and is zero otherwise.
The following lemma \cite{BNNaT}, which is an immediate consequence of the infinitesimal braid relations,
will play a key role.
\begin{lem}[Naturality of cabling]\label{natCabling}
For any positive integers $a$, $b$, $c$, and any $x\in \widehat{\mathcal{A}}_{b}$, $y\in \widehat{\mathcal{A}}_{a+1+c}$, one has the equality in figure \ref{naturalCabling}.
\end{lem}
We will also need the coproduct
map
 $$\square : \widehat{\mathcal{A}}_n\rightarrow \widehat{\mathcal{A}}_n \widehat{\otimes} \widehat{\mathcal{A}}_n$$
defined so that it is continuous and each generator $t_{ij}$ is primitive. An element
$\Psi\in \widehat{\mathcal{A}}_n$
is group-like if and only if $\square \Psi=\Psi\otimes \Psi$ ; this implies invertibility.
 We denote by $\mathcal{T}_n$ the Lie algebra of primitive elements in $\mathcal{A}_n$, and by
$\widehat{\mathcal{T}}_n$ its closure in $\widehat{\mathcal{A}}_n$. The Lie algebra  $\mathcal{T}_n$
is known as the Lie algebra of infinitesimal braids. 

\begin{dfn}
An { associator} is a group-like
element $\Phi\in \widehat{\mathcal{A}}_3$
satisfying
\begin{itemize}
\item the pentagon identity
$$(\un\otimes \Phi).d_2\Phi.(\Phi\otimes \un)=d_3\Phi.d_1\Phi$$
\item the hexagon identities
$$d_1 \exp\left(\pm \frac{1}{2}t_{12}\right)={}^{312}\! \Phi \exp\left(\pm \frac{1}{2}t_{13}\right)
{}^{132}\!\!\left(\Phi^{-1}\right)\exp\left(\pm \frac{1}{2}t_{23}\right)\Phi 
$$
(The notation $ijk$
is used for the permutation $(1,2,3)\mapsto (i,j,k)$.)
\item $\Phi$ is non-degenerate: $s_1\Phi=s_2\Phi=s_3\Phi=\mathbf{1}_{(\bullet\bullet)}$;
\end{itemize}
 \end{dfn}
\begin{rems}
1. The group-like element $\Phi$ can be written
$\Phi=\exp(\phi)\ \ \text{, with }\phi\in \widehat{\mathcal{T}}_3\ .$
From the non degeneracy condition, we get that $\phi$ belongs to $[\widehat{\mathcal{T}}_3,\widehat{\mathcal{T}}_3]$.\\
2. The hexagon identity with negative sign could be replaced by
${}^{321}\Phi=\Phi^{-1}\ .$
3. The non degeneracy condition can be deduced from group-like and pentagone properties \cite[Proposition 1]{REPTHEO}
\end{rems}

 The notion of associator is due to Drinfeld \cite{DRINFELD}, who also obtained an associator with complex
 coefficients using monodromy of KZ system, showed that associators with rational coefficients exist,
 and emphasized the role of associators in constructing monoidal functors from braids to the universal envelopping algebras of infinitesimal braids.
  This functor was extended to the tangle category by Kontsevich
  \cite{KONTSEVICH} and further developped by Bar-Natan, Le-Murakami
  and others.
  Considering parenthesized braids and tangles converts Drinfeld-Kontsevich functor into a strictly monoidal functor.
 Drinfeld-Kontsevich functor
 on the whole parenthesized tangle category is constructed in \cite{LM1,BNNaT}. 
 We  consider as a target category the (strict) monoidal category $\widehat{\mathfrak{B}}$
whose objects are integers and morphisms are defined by
$$\End_{\widehat{\mathfrak{B}}}(n)=\widehat{\mathfrak{B}}_n\ ,$$
$$\mathrm{Hom}_{\widehat{\mathfrak{B}}}(n,m)=\{0\}\ \  \text{ for }  n\neq m\ .$$

\vspace{5pt}
\noindent{\bf Cabling associators and braidings.} 
 For $p\in \mathrm{Obj}(\mathbf{PaB})$, we denote by $|p|$ the number of points of the parenthesization $p$.
  Let $\Phi$ be an associator and let $p,q,r \in \mathrm{Obj}(\mathbf{PaB})$, we define the cabled associator
 $\Phi_{p,q,r}$ by
$$\Phi_{p,q,r}=(d_1^{|p|}\circ d_2^{|q|}\circ d_3^{|r|}) \Phi\ .$$
 The identity braid, viewed as an element
in $\mathrm{Hom}_{\mathbf{PaB}}\left((pq)r,p(qr)\right)$, is denoted by $a_{p,q,r}$.

Let $\mathrm{R}\in \widehat{\mathcal{B}}_2$ be the element depicted below. Here a product from  right to left
is depicted from the bottom.

\centerline{
{\Png{2.5cm}{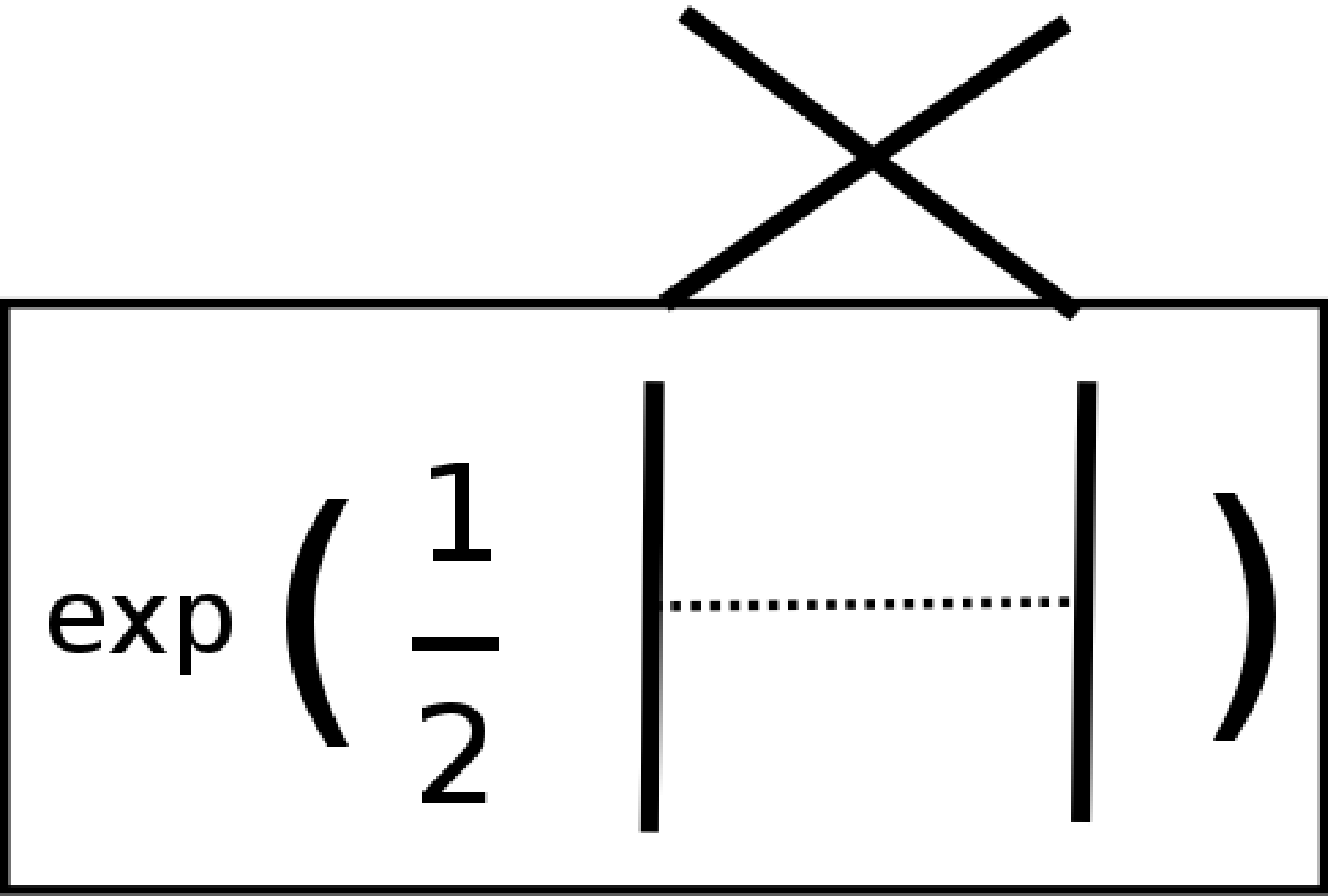}}}

For parenthesizations $p$ and $q$, $\mathrm{R}_{p,q}$ denote
$d_1^{|p|}d_2^{|q|}\mathrm{R}$, and $c_{p,q}\in \mathrm{Hom}_{\mathbf{PaB}}\left(pq,qp\right)$ denote the element depicted below.

\centerline{
{\Png{3cm}{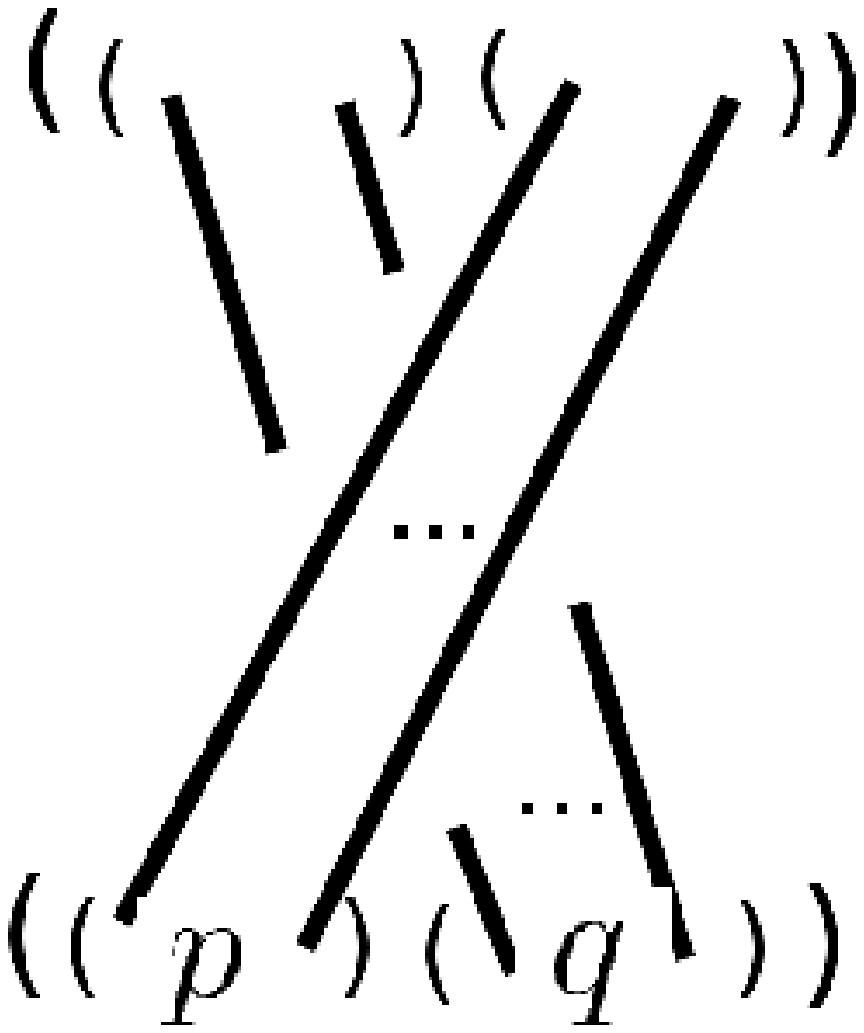}}}

From $d_i^2=d_{i+1}d_i$, we get that $(p,q,r)\mapsto \Phi_{p,q,r}$ and $(p,q)\mapsto \mathrm{R}_{p,q}$ commutes with the cabling operation. We quote that $\Phi_{p,q,r}$ and $\mathrm{R}_{p,q}$ only depend on the integers $|p|$, $|q|$, $|r|$.

\vspace{5pt}
\begin{theo}
If $\Phi$ is an associator, then there exists a unique strictly monoidal functor
$$Z :\mathbf{PaB}\rightarrow \widehat{\mathfrak{B}}$$
such that, for a parenthesization $p$, $Z(p)=|p|$, and for all parenthesizations $p$, $q$, $r$
$$Z(c_{p,q})=\mathrm{R}_{p,q}\ \ ,\ \ Z(a_{p,q,r})=\Phi_{p,q,r}\ .$$
%
\end{theo}
From the  definition of the functor $Z$ in the above theorem, we  
obtain the following result, first proved by Le-Murakami \cite{LM2}.
\begin{ccor}
The functor $Z$ commutes with cabling operations.
\end{ccor}
\begin{proof}
The elements    $c_{p,q}$ and  $a_{p,q,r}$ generate $\mathbf{PaB}$ as a monoidal category; unicity of the functor $Z$ follows.

As already said, the proof of existence can be found in \cite{LM1,BNNaT} and rests on Mac-Lane coherence theorem
for braided monoidal categories \cite[Section 2]{JS}. We like to give here a rather self contained argument.
 From the pentagon identity, we get that for any two parenthesizations $P$, $Q$
with the same length $|P|=|Q|=n$ we have a canonical invertible element $\Phi_P^Q \in\hat{\mathcal{A}}_n$
obtained by composing associators corresponding to a decomposition of the identity braid,
considered as an element in $\mathrm{Hom}_{\mathbf{PaB}}(P,Q)$, into tensor products and composition of
elements $a_{p,q,r}^{\pm1}$.

Now consider a parenthesized braid $\tau\in\mathrm{Hom}_{\mathbf{PaB}}(P,Q)$ , and represent it as a diagram
corresponding to a word $w$ in Artin generators $w=\sigma_{i_k}^{\pm1}\dots \sigma_{i_1}^{\pm1}$.
 For each integer $n$ we denote by $p(n)$ the parenthesization {\em from the left}, i.e.
$p(n+1)=(p(n)\bullet )$, and $p_i=p_i(n)=d_i(p(n-1))$.
  We decompose $\tau$ as follows
$$\tau=a^Q_{p_{i_k}} \sigma_{i_k}^{\pm1}a_{p_{i_{k-1}}}^{p_{i_k}}\dots \sigma_{i_1}^{\pm1}a_P^{p_{i_1}}\ .$$
For $0<i<n$, let $\mathrm{R}_i=s_i\exp(\frac{1}{2}t_{i,i+1})\in{\widehat{\mathcal{B}}_n}$
($\mathrm{R}$ is inserted along strands $i$, $i+1$).
 Here $s_i\in\mathfrak{S}_n$ is the transposition $(i,i+1)$.

 We associate to the representative of $\tau$ written above
$$Z(\tau,w)=\Phi^Q_{p_{i_k}} \mathrm{R}_{i_k}^{\pm1}\Phi_{p_{i_{k-1}}}^{p_{i_k}}\dots \mathrm{R}_{i_1}^{\pm1}\Phi_P^{p_{i_1}}\ .$$
We will show that the assignement $t\mapsto Z(t)=Z(t,w)$ constructs a well defined functor
$Z$ with the required properties. This follows from the statements below whose proof is left to the reader.
 This is just playing with the definition, hexagon and pentagon relations, together with the naturality property
of the cabling operation on chord diagrams (lemma \ref{natCabling}).
\begin{enumerate}
\item Functoriality: $Z(\tau \tau',ww')=Z(\tau,w)Z(\tau',w')$, for any $\tau$, $\tau'$ for which $\tau\tau'$ makes sense.
\item Braid relations: if $w$ and $w'$ are equivalent as braids, then we have $Z(\tau,w)=Z(\tau,w')$.
Hence we get that $Z$ is a well defined functor on $\mathbf{PaB}$.\\
It is sufficient to show that for $0<i+1<j<n$,
$$Z(\tau,w_1\sigma_i\sigma_jw_2)=Z(\tau,w_1\sigma_j\sigma_iw_2) \ ;$$
and for $0<i<n$,
$$Z(\tau,w_1\sigma_i\sigma_{i+1}\sigma_iw_2)=Z(\tau,w_1\sigma_{i+1}\sigma_i\sigma_{i+1}w_2) \ .$$
\item $Z$ is strictly monoidal, i.e. $Z$ is compatible with juxtaposition.
\item $Z(c_{p,q})=\mathrm{R}_{p,q}$, and $Z(a_{p,q,r})=\Phi_{p,q,r}$.
\end{enumerate}
\end{proof}
\subsection{Commutation property}
We already considered the cabling operator which duplicates one strand.
 We want now to duplicate all strands. We denote by $\Delta$
the functor $\mathbf{PaB}\rightarrow \mathbf{PaB}$ which duplicate all strands
(and put a parenthesis around each pair). More generally, we denote by
$\Delta_r$ the functor $\mathbf{PaB}\rightarrow \mathbf{PaB}$ which
replace each strand by $r$ parallel copies (for each $r$-uple put parenthesis from the left).
The functors for chord diagrams similar to $\Delta$ and $\Delta_{r}$ are denoted respectively
by $\partial$ and $\partial_{r}$. We saw that the functor $Z$ commutes with the cabling operation. This shows the following proposition.
\begin{prop}\label{squareFunctor}
The following square of functors is commutative.
$$
\xymatrix{
\mathbf{PaB} \ar[d]^{\Delta_{r}} \ar[r]^{Z} & \WB\ar[d]^{\partial_{r}} \\
\mathbf{PaB} \ar[r]^{Z} &  \WB
}
$$

\end{prop}

  For each parenthesization $p$ with $|p|=n$,
 by extending linearly the map, we have an algebra morphism
$x\mapsto Z(x)$:
$$Z_p: \k B_n\approx \k End_{\mathbf{PaB}}(p)\rightarrow \widehat{\mathfrak{B}}_n\ .$$
  If $q$ is another parenthesization with $|q|=n$, then $Z_q$ is obtained by intertwining with
 $Z(\text{identity braid})$, where identity braid is considered as an element in $\mathrm{Hom}_{\mathbf{PaB}}(p,q)$.
 
\begin{dfn}
Let $\sigma_i$, $i \in [1,n-1]$, be the Artin
generators of $B_n$.
A morphism $\varphi : \k B_n \to \widehat{\B}_n$ for some
$n$ is said to satisfy property (*) if, for all $i \in [1,n-1]$,
$$
(*) \ \ \ \ \varphi(\sigma_i) \mbox{ is conjugated to }
\mathrm{R}_i=s_i\exp(\frac{1}{2}t_{i,i+1})\mbox{ by
some }\psi_i \in \exp [\widehat{\mathcal{T}_n},\widehat{\mathcal{T}_n}]
$$
Here $s_i\in\mathfrak{S}_n$ is the transposition $(i,i+1)$.
\end{dfn}
This implies that
$\varphi(\sigma_i)$ equals $s_i(1+\frac{1}{2} t_{i,i+1})$
plus higher terms.

  Recall that for an integer $n$, we denote  by $p(n)$ the object in $\mathbf{PaB}$ equal to the left
parenthesization $(\dots((\bullet\bullet)\bullet)\dots \bullet)$.
 Observe that the object $\Delta_r(p(n))$ is not equal to the left
parenthesization 
$p(rn)$, but rather: $\Delta_rp(n)=(\dots((p(r))p(r))\dots p(r))$.
 We denote by $\Delta_{r,n}$ the linear extension of
$\Delta_{r}p(n)$:
$$\Delta_{r,n}: \k B_n\approx \k \mathrm{End}_{\mathbf{PaB}}(p(n))
\rightarrow \k\mathrm{End}_{\mathbf{PaB}}(\Delta_rp(n))\approx \k B_{rn}\ .$$

The morphism for chord diagrams similar to $\Delta_{r,n}$ is denoted
by $\partial_{r,n}$.

\begin{prop} \label{commutation} For all $n\geq 1$, $r \geq 2$, the
morphisms $$\varphi_n=Z_{n(p)}: \k B_n \approx \k \mathrm{End}_{\mathbf{PaB}}(p(n))
\to \WB_n$$
and  $$\varphi_{r,n}=Z_{\Delta_r(n(p))} : \k B_{nr} \approx \k \mathrm{End}_{\mathbf{PaB}}(  \Delta_r(n(p)) )\to \WB_{nr}$$
satisfy (*) and  the following diagram commutes :
$$
\xymatrix{
\k B_n \ar[d]^{\Delta_{r,n}} \ar[r]^{\varphi_n} & \WB_n \ar[d]^{\partial_{r,n}} \\
\k B_{nr} \ar[r]^{\varphi_{r,n}} &  \WB_{nr}
}
$$
\end{prop}
\begin{proof}
We have property (*), with $\psi_i=\Phi_{n(p)}^{p_i}$.
 Commutativity of the diagram follows from proposition \ref{squareFunctor}.
\end{proof}

\section{Representations of the braid groups from an infinitesimal view
point}

For the decomposition of representations of the braid groups $B_n$
with generic parameters it is much more easier, if possible, to deal
with the representations of $\B_n$ from which they arise. Indeed,
there is a dictionnary between properties of representations
of $\B_n$ and the representations of $B_n$ obtained
through a morphism $\varphi$ such as the ones discussed above.

This dictionnary is discussed at length in \cite{REPTHEO}. In the
following sections we show how to get the representations
involved here from `` infinitesimal representations ''.

\subsection{General facts}

Let $\K = \k((h))$ and $\L = \k(q)$. We view $\L$ as embedded in
$\K$ by $q \mapsto \exp( h/2)$. Let $n \geq 2$ and
$\varphi : \k B_n \to \WB_n$ satisfying (*).
To every representation $\rho : \B_n \to\End(\k^n)$ one associates
a representation $\tilde{\rho} : \WB_n \to \End(\K^n)$
defined by
$$
\left\lbrace
\begin{array}{lclr}
\tilde{\rho}(t_{ij}) & = & h \rho(t_{ij}) \\
\tilde{\rho}(s) & = & \rho(s) & \mbox{ if } s \in \SS_n.
\end{array} \right.
$$
and a representation $\wp(\rho) = \tilde{\rho}\circ \varphi$.
of $B_n$. The condition (*) implies that
$\rho$ is (absolutely) irreducible iff $\wp(\rho)$ is
(absolutely) irreductible (cf. \cite{REPTHEO}, proposition 7 and proposition 8).

In general, for any representation $\rho$ of $\B_n$ and $R = \wp(\rho)$
with $\varphi$
satisfying (*), one has $R(\sigma_i) \equiv \rho(s_i)$ modulo $h$. It
follows that $\Ker R \subset P_n$ as soon as the restriction of $\rho$
to $\mathfrak{S}_n$ is faithful. Let us assume that this is the case. Recall
that $P_n$ has standard generators
$\xi_{ij}$ with $1 \leq i < j \leq n$ whose images generate $P_n/(P_n,P_n)$.
Under condition (*) one has $R(\xi_{ij}) = 1 +  h \rho(t_{ij})$
modulo $h^2$, hence
$\Ker R \subset (P_n,P_n)$ as soon as the endomorphisms
$\rho(t_{ij})$ are linearly independant.

\subsubsection{Framing}

For any given representation $R : B_n \to GL_N(\K)$ and $a \in \K^{\times}$
one may define another representation, denoted $a R$ and defined
by $(a R)(\sigma_i) = a R(\sigma_i)$ for all $i \in [1,n-1]$.
It is clear that, if $\Ker R \subset (B_n,B_n)$,
then $\Ker  R \subset \Ker a R$.

At the infinitesimal level, to any $\rho : \B_n \to \End(\k^N)$
and $v \in \K$ one can associate the representation $v + \rho$
defined by
$$
\left\lbrace
\begin{array}{lclr}
(v+\rho)(t_{ij}) & = &  \rho(t_{ij}) + v \Id \\
(v+\rho)(s) & = & \rho(s) & \mbox{ if } s \in \SS_n.
\end{array} \right.
$$

It is easily checked that, under condition (*), one has
$$
\wp(v+ \rho) = \exp(hv) \wp(\rho).
$$

Let us assume that $R = \wp(\rho)$
with $\varphi$ satisfying (*). Then $\det R(x)$ is a formal series
in $h$ with constant term 1. Such series form a group with unique roots.
When $\Ker R \subset P_n$,
a sufficient condition for the converse inclusion $\Ker  a R 
\subset \Ker  R$ to hold is then that
$\det R(\sigma_1)^2 \neq a^{- 2\dim R}$. Indeed, if $l = l(x)$ denotes
the length of $x \in B_n$ with respect to the set of Artin generators,
$x \in \Ker a R$ means that $R(x) = a^{-l}$, hence $\det R(x) = a^{-l \dim R}$.
If $l =0$, then $x \in \Ker R$ and we are done. Otherwise,
since the Artin generators are conjugated one to the other,
then $\det R(x) = \det(R(\sigma_1))^l$, thus $\det(R(\sigma_1^2))^l = (a^{-2 \dim R})^l$
hence $\det R(\sigma_1)^2 = a^{-2 \dim R}$, which has been ruled out.

\subsubsection{Twisting}

For any given representation $R : B_n \to GL_N(\L)$ and $r \in \Z \setminus
\{ 0 \}$
one may also define another representation, denoted $R^{q^r}$
and defined by twisting $R$ by the field morphism $q \mapsto q^r$.
It is clear that $R$ and $R^{q^r}$ have the same kernel.

In a similar manner, to any $b \in \k^*$ and any representation $\rho$ of
$\B_n$ one can associate the
representation $b \rho$ defined by
$$
\left\lbrace
\begin{array}{lclr}
(b \rho)(t_{ij}) & = &  b \rho(t_{ij})  \\
(b \rho)(s) & = & \rho(s) & \mbox{ if } s \in \SS_n.
\end{array} \right.
$$

In the situations described below, these two operations are intimately
related.

\subsection{The Burau representation}

Let $H_n(q) = \L B_n / I$, where $I$ is the ideal generated
by the elements $(\sigma_i - q)(\sigma_i + q^{-1})$ for
$i \in [1,n]$. Note that, since the Artin generators are
conjugated one to the other, this is the same as the ideal
generated by $(\sigma_1 - q)(\sigma_1 + q^{-1})$. This
algebra is the (generic) Iwahori-Hecke algebra of type $A$.
This is a well-known finite-dimensional algebra, isomorphic
to the group algebra over $\L$ of the symmetric group
on $n$ letters. Its representations are explicitely described
in \cite{HOEFSMIT,WENZL}.

Let $V$ be a $n$-dimensional $\L$-vector space with basis
$e_1,\dots,e_n$. The Burau representation of the braid group
$\Bur : B_n \to GL(V)$ is defined in matrix block-diagonal
form on this basis by
$$
\Bur(\sigma_i) = 
q I_{k-1} \oplus \left( \begin{array}{cc} q-q^{-1} & q \\
q^{-1} & 0 \end{array} \right) \oplus q I_{n-k+1}.
$$
It is easily checked that $(\Bur(\sigma_1) - q)(\Bur(\sigma_1) + q^{-1}) =
0$. From the classical representation theory of the Iwahori-Hecke
algebra, one checks that the Burau representation is
characterized among its representations by the following properties
\begin{enumerate}
\item It is $n$-dimensional, with a $(n-1)$-dimensional
irreducible subspace
\item $-q^{-1}$ has multiplicity $1$ in the spectrum of $\sigma_1$.
\end{enumerate}
In particular, the image of $\sigma_1$ has determinant $-q^{n-2}$.
The Burau representation $\Bur$ can be deduced from the
representation $\bur$ of $\B_n$ defined by
$$
\left\lbrace
\begin{array}{rclr}
\bur(t_{ij}).e_i & = & e_j \\
\bur(t_{ij}).e_k & = & e_k & \mbox{ if } k \not\in \{i, j\} \\
\bur(s).e_i & = & e_{s(i)} & \mbox{ if } s \in \SS_n
\end{array} \right.
$$
Indeed, given any $\varphi : \k B_n \to \WB_n$ satisfying (*),
the representation $\wp(\bur)$ extended to $\L$ factorizes
through $H_n(q)$ because the eigenvalues of $s_1$ are
$1,-1$ and $\wp(\bur)(\sigma_1)$ is conjugated to
$\bur(s_1) \exp( h \bur(s_1)/2)$. Moreover, $-1$ has
multiplicity $1$ in the spectrum of $\bur(s_1)$, hence $-q^{-1}$
has multiplicity $1$ in the spectrum of $\wp(\bur)(\sigma_1)$.
Finally, the kernel of the linear map $\alpha : V \to L$
defined by $\alpha(e_i) =1$ is stable under $\bur$, and
it is a standard fact from the representation theory
of the symmetric group that it is irreducible under
the action of $\SS_n$. It follows that $\wp(\bur)$
admits an irreducible $(n-1)$-dimensional subspace, hence
the following proposition.
\begin{prop} For all $\varphi : \k B_n \to \WB_n$ satisfying
(*), the Burau representation is isomorphic to $\wp(\bur)$.
\end{prop}
Moreover, by the same characterization, one gets that, once extended
over $\K$, the twisted representation $\Bur^{q^r}$ is isomorphic
to $\wp(r \bur)$ for all $r \in \Z\setminus \{0 \}$.

It is readily checked that the restriction of $\Bur$ to the
center is faithful for $n \geq 3$. For $n \geq 5$ it
is known to be unfaithful. Indeed, Bigelow found in \cite{BIGELOW},
improving results of Moody \cite{MOODY}, Long and Paton \cite{LONGPATON}, that the element
$\beta = (\psi_2 \psi_1^{-1} \sigma_4 \psi_1 \psi_2^{-1},\delta_5)$
is non trivial and lies in the kernel of $\Bur$ for $n \geq 5$,
where
$$
\left\lbrace \begin{array}{lcl}
\psi_1 &=& \ss_3^{-1} \ss_2 \ss_1^2 \ss_2 \ss_4^3 \ss_3 \ss_2\\
\psi_2 & = & \ss_4^{-1} \ss_3 \ss_2 \ss_1^{-2} \ss_2 \ss_1^2 \ss_2^2 \ss_1 \ss_4^5\\
\delta_5 & = & \ss_4 \ss_3 \ss_2 \ss_1^2 \ss_2 \ss_3 \ss_4. 
\end{array}
\right.
$$
At the present time, very few things are known about this kernel.
However, from the basic observations above, the faithfulness of the
permutation representation of the symmetric group implies that
$\Ker \Bur \subset P_n$, and the fact that the $\rho(t_{ij})$
are linearly independent implies $\Ker \Bur \subset (P_n,P_n)$.

\subsection{The extended permutation representation}

The \emph{extended symmetric group} (from the french ``groupe sym\'etrique
\'etendu'') can be defined as $\ST_n = B_n / (P_n,P_n)$ and has been
introduced by J. Tits in \cite{TITS}. If $\rho$ is a representation
of $\B_n$ and $\varphi : \K B_n \to \WB_n$ satisfies
(*), then $\wp(\rho)$ factorizes through $\ST_n$ if and only if
$\rho([t_{ij},t_{kl}]) = 0$ for all $i,j,k,l$, which is equivalent
to saying that $\rho([t_{12},t_{23}]) = 0$ (see \cite{REPTHEO} lemma 5).
Moreover, in that case, for all $i \in [1,n]$,
$$
\wp(\rho)(\sigma_i) = \rho(s_i) \exp( h \rho(t_{i,i+1})/2).
$$
The natural representation $\Sym$ of $\ST_n$ over $\L$
can be defined over a $n$-dimensional vector space $\L^n$ in
block-diagonal form by
$$
\Sym(\sigma_k) = I_{k-1} \oplus \left( \begin{array}{cc} 0 & q \\
q & 0 \end{array} \right) \oplus I_{n-k+1}.
$$
It is readily checked that, if $\varphi$ satisfies (*) then
$\Sym = \wp(\sym)$ where the restriction of $\sym$ to
$\SS_n$ is defined in
block-diagonal form by
$$
\begin{array}{lcl}
\sym(s_k) & = & I_{k-1} \oplus \left( \begin{array}{cc} 0 & 1 \\
1 & 0 \end{array} \right) \oplus I_{n-k+1} \\
\sym(t_{k,k+1}) & = & 0_{k-1} \oplus \left( \begin{array}{cc} 1 & 0 \\
0 & 1 \end{array} \right) \oplus 0_{n-k+1}.
\end{array}
$$
More generally, $\wp(r \sym) = \Sym^{q^r}$.
It is easily checked that $\sym$, thus $\Sym$, is
(absolutely) irreducible (cf. \cite{THESE} {\bf II} 2.3.1 lemme 8).
Again by faithfulness of the permutation representation of the symmetric
group one gets that $\Ker \Sym \subset P_n$. Moreover by definition
of the extended symmetric group one has $(P_n,P_n) \subset \Ker \Sym$,
hence $\Ker \Bur \subset \Ker \Sym$.

\section{Cabling the Burau representation}

\subsection{Infinitesimal result}

Let us consider the infinitesimal Burau representation of $\B_{nr}$.
We can index the basis elements $e_s$ as $e_i^j$ for $1 \leq i \leq r$
and $1 \leq j \leq n$ with $s = r(j-1)+i$. Using $\partial_{n,r} :
\B_n \to \B_{nr}$ one gets a representation $\bur \circ \partial_{n,r}$
defined by $s.e_i^j = e_i^{s(j)}$ for $s \in \SS_n$ and
$$
\left\lbrace \begin{array}{lll}
t_{ij} . e_s^i & = \sum_{t \in [1,r]} e_t^j  +  r(r-1) e_s^i& \\
t_{ij} . e_s^j & = \sum_{t \in [1,r]} e_t^i +  r(r-1) e_s^j & \\
t_{ij} . e_s^k & =  r^2 e_s^k & \mbox{ if } k \not\in \{i,j \}. \\
\end{array}
\right.
$$
For all $1 \leq i \leq n$, Let us introduce
$u_i = \sum_{s \in [1,r]} e_s^i$. The subspace
$U$ generated by $u_1,\dots,u_n$ has dimension $n$. From the above formulas
one gets
$$
\left\lbrace \begin{array}{lll}
t_{ij} . u_i & =  r u_j +  r(r-1) u_i & \\
t_{ij} . u_j & =  r u_i + r(r-1) u_j  & \\
t_{ij} . u_k & =  r^2 u_k & \mbox{ if } k \not\in \{i,j \}. \\
\end{array}
\right.
$$
This means that, on $U$, $t_{ij}$ acts in the same way that
$
r  (i \ j) +
 r(r-1)$, hence the action of $\WB_n$ is isomorphic to $r\bur + r(r-1)$.
Using $\varphi$ one thus gets a representation of $B_n$ isomorphic to
$r(r-1)\Bur^{q^r}$.

Let $E_i$ for $1 \leq  i \leq n$ denote the subspace spanned by
$e_s^i$ pour $s \in [1,r]$, and $\alpha_i$ be the linear form on
$E_i^*$ defined by
$\alpha_i(e_s^i) = 1$. Let $K_i = \Ker \alpha_i \subset E_i$,
which is spanned by the $e_s^i-e_t^i$, and let $K$ be the (direct) sum
of these subspaces.

If $x \in K_i$ or $x \in K_j$, one has
$t_{ij}.x =  r(r-1) x$ and, if $x \in K_k$ for $k \not\in \{ i, j \}$,
then $t_{ij}.x = 
r^2$. It follows that $K$ is stable under the action of
$\B_n$ and that the action of $\mathcal{A}_n$ on $K$ is commutative.
More precisely, on the subspaces $F_{s,t}$ spanned
by the $v^i_{s,t} = e^i_s - e^i_t$ for $i \in [1,n]$,
when $s \neq t$, the action of $\B_n$ is isomorphic to $r^2-r
\sym$. It follows that we have the following result.
\begin{prop} The representation $\bur \circ \partial_{n,r}$
is isomorphic to the direct sum of $r \bur + r(r-1)$ and
$r-1$ copies of $r^2-r
\sym$.
\end{prop}

\subsection{Global result}

Since there exist morphisms $\varphi$ and $\varphi_{n,r}$ satisfying
the conclusions of proposition \ref{commutation}, the above
proposition implies our main theorem.
\begin{theo}
The representation $\Bur \circ \Delta_{n,r}$
is isomorphic to the direct sum of $ q^{r(r-1)}\Bur^{q^r}$ and
$r-1$ copies of $q^{r^2}
\Sym^{q^{-r}}$.
\end{theo}

As far as faithfulness questions are concerned, this theorem
implies that the kernel of $\Bur \circ \Delta_{n,r}$ in $B_n$ is the
intersection of those of $q^{r(r-1)}\Bur$ and $q^{r^2}\Sym$. 
Because $q^{2n-4} \neq q^{-2r(r-1)n}$ for all values of $r \geq 1$,
and $\Ker \Bur \subset (B_n,B_n) \cap P_n$, then $\Ker q^{r(r-1)}\Bur
= \Ker \Bur$. Moreover $\Ker \Bur \subset (P_n,P_n) \subset
\Ker  q^{r^2} \Sym$. It follows that cabling does not
improve faithfulness of the Burau representation.

\begin{ccor} For all $n \geq 3$ and $r \geq 1$, the kernel of $\Bur \circ \Delta_{n,r}$ coincides with the
kernel of $\Bur$. In particular
$\Bur \circ \Delta_{n,r}$ is not faithful for $n \geq 5$.
\end{ccor}


\begin{thebibliography}{9999}
\bibitem[BN1]{BNNaT} Bar-Natan, Non-associative tangles, in ``Geometric topology (Athens, GA, 1993)'' 139-183 AMS/IP Stud. Adv. Math., 2.1, Amer. Math. Soc., Providence, RI (1997).
\bibitem[BN2]{BNGT} D. Bar Natan, {\it On associators and the Grothendieck-Teichmuller group},
Sel. math., New ser. 4, 183--212 (1998).
\bibitem[Bi]{BIGELOW} S. Bigelow, {\it The Burau representation is not faithful for $n=5$},
Geom. Topol. 3  397--404 (1999).
\bibitem[D]{DRINFELD} V. G. Drinfeld, {\it On quasitriangular quasi-Hopf algebras and a group closely connected with $\mathrm{Gal}(\overline{\Q}/\Q)$},
Leningrad Math. J. 2, 829--860 (1991).
\bibitem[Ho]{HOEFSMIT} P.N. Hoefsmit, {\it Representations of Hecke algebras of finite groups with BN-pairs of classical type},
Ph.D. Thesis, University of British Columbia, 1974.
\bibitem[JS]{JS} A. Joyal and R. Street, {\it Braided tensor categories}, Advances in Math. 102, 20--78 (1993).
\bibitem[K]{KONTSEVICH} M. Kontsevich,
{\it Vassiliev's knot invariants.}
 Adv. Sov. Math. 16(2), 137-150 (1993).
\bibitem[LM1]{LM1} T. T. Q. Le and J. Murakami, {\it The universal Vassiliev-Kontsevich invariant for framed oriented links}, Compositio Math. 102, 42--64 (1996). 
\bibitem[LM2]{LM2} T. T. Q. Le and J. Murakami, {\it Parallel version of the universal Vassiliev-Kontsevich invariant}, J. Pure and Applied Algebra 121, 271--291 (1997).
\bibitem[L]{LESCOP} C. Lescop, {\it Introduction to the Kontsevich integral of framed tangles}, Summer School lectures, Grenoble 1999, \verb+www-fourier.ujf-grenoble.fr/ECOLETE+
\bibitem[LP]{LONGPATON} D.D. Long, M. Paton, {\it The Burau representation is not faithful for $n\geq 6$}, Topology {\bf 32}, no. 2, 439--447  (1993).
\bibitem[Mo]{MOODY} J.A. Moody, {\it The Burau representation of the braid group $B\sb n$ is unfaithful for large $n$}, Bull. Amer. Math. Soc. (N.S.) {\bf 25}, no. 2, 379--384 (1991).
\bibitem[Ma1]{THESE} I. Marin, {\it Repr\'esentations lin\'eaires
des tresses infinit\'esimales}, Ph. D. Thesis, Orsay, 2001.
\bibitem[Ma2]{REPTHEO} I. Marin, {\it On the representation theory
of braid groups}, preprint \verb+arXiv:math.RT/0502118 (v3)+.
\bibitem[Ti]{TITS} J. Tits, {\it Normalisateurs de tores I. Groupes de Coxeter \'etendus}
J. Algebra 4, 96--116 (1966).
\bibitem[Wz]{WENZL} H. Wenzl, {\it Hecke algebras of type $A\sb n$ and subfactors},
Invent. Math. 92  349--383 (1988).

\end{thebibliography}
\end{document}